\documentclass[12pt]{article}
\usepackage{pstricks,subfigure,graphicx}

\usepackage{amsmath,amssymb}
\usepackage{tikz}
\usepackage{color}

\newtheorem{theorem}{Theorem}[section]
\newtheorem{lemma}[theorem]{Lemma}
\newtheorem{corollary}[theorem]{Corollary}

\newtheorem{problem}[theorem]{Problem}

\newcommand{\proof}{\noindent{\bf Proof.\ }}
\newcommand{\qed}{\hfill $\square$ \bigskip}
\newcommand{\GG}{G\overline{G}}

\def\cp{\,\square\,}
\DeclareMathOperator {\gp} {gp}

\DeclareMathOperator {\diam} {diam}
\DeclareMathOperator {\rad} {rad}
\DeclareMathOperator {\ecc} {ecc}
\let\deg\relax
\DeclareMathOperator {\deg} {deg}

\begin{document}

\title{On the general position number of complementary prisms}

\author{
    Neethu P. K.$^{a}$
    \and
    	Ullas Chandran S. V. $^{a}$
    	\and
    		Manoj Changat$^{b}$
    		\and
    Sandi Klav\v{z}ar$^{c,d,e}$ 
    	}
	
\date{\today}

\maketitle

\begin{center}
	$^a$ Department of Mathematics, Mahatma Gandhi College, Kesavadasapuram,  Thiruvananthapuram-695004, Kerala, India \\
	{\tt p.kneethu.pk@gmail.com}, {\tt svuc.math@gmail.com} 
	\medskip
	
	$^b$ Department of Futures Studies, University of Kerala Thiruvananthapuram-695034, Kerala, India \\  
	{\tt mchangat@gmail.com} 
	\medskip
	
	$^c$ Faculty of Mathematics and Physics, University of Ljubljana, Slovenia\\
	{\tt sandi.klavzar@fmf.uni-lj.si}
	\medskip
	
	$^d$ Faculty of Natural Sciences and Mathematics, University of Maribor, Slovenia
	\medskip
	
	$^e$ Institute of Mathematics, Physics and Mechanics, Ljubljana, Slovenia
\end{center}

\begin{abstract}
The general position number  ${\rm gp}(G)$ of a graph $G$ is the cardinality of a largest set of vertices $S$ such that no element of $S$ lies on a geodesic between two other elements of $S$. The complementary prism $G\overline{G}$ of $G$ is the graph formed from the disjoint union of $G$ and its complement $\overline{G}$ by adding the edges of a perfect matching between them. It is proved that ${\rm gp}(\GG)\le n(G) + 1$ if $G$ is connected and ${\rm gp}(\GG)\le n(G)$ if $G$ is disconnected. Graphs $G$ for which ${\rm gp}(\GG) = n(G) + 1$ holds, provided that both $G$ and $\overline{G}$ are connected, are characterized. A sharp lower bound on ${\rm gp}(\GG)$ is proved. If $G$ is a connected bipartite graph or a split graph then ${\rm gp}(\GG)\in \{n(G), n(G)+1\}$. Connected bipartite graphs and block graphs for which ${\rm gp}(G\overline{G})=n(G)+1$ holds are characterized. A family of block graphs is constructed in which the ${\rm gp}$-number of their complementary prisms is arbitrary smaller than their order.
\end{abstract}


\medskip\noindent
{\bf AMS Subj.\ Class:} 05C12; 05C69

\section{Introduction}
\label{sec:intro}

The general position problem in graphs was introduced in~\cite{manuel-2018a} as a graph theory variant of the classical, century old Dudeney's no-three-in-line problem~\cite{dudeney-1917} and the general position subset selection problem from discrete geometry~\cite{froese-2017, ku-2018, misiak-2016, payne-2013, skotnica-2019}. A set $S$ of vertices in a graph $G$ is a {\em general position set}  if no element of $S$ lies on a geodesic between any two other elements of $S$. A largest general position set is called a {\em $\gp$-set} and its size is the {\em general position number} ({\em $\gp$-number} for short) \ ${\rm gp}(G)$ of $G$. The same concept was in use two years earlier in~\cite{ullas-2016} under the name {\em geodetic irredundant  sets}.  

Let us briefly recall the progress on the general position problem so far. In~\cite{manuel-2018a}, general upper and lower bounds on the $\gp$-number were proved as well as NP-completeness of the problem for arbitrary graphs.  The ${\rm gp}$-number of a large class of subgraphs of the infinite grid graph and of some other classes were obtained in~\cite{manuel-2018b}. The paper~\cite{bijo-2019} gives a characterization of general position sets which is then applied in  determining the ${\rm gp}$-number of graphs of diameter $2$, cographs, graphs with at least one universal vertex, bipartite graphs and their complements. Subsequently, the $\gp$-number of complements of trees, of grids, and of hypercubes were deduced  in~\cite{bijo-2019}. In~\cite{ghorbani-2019}, a sharp lower bound on the ${\rm gp}$-number of Cartesian products is proved, and the ${\rm gp}$-number for different graph operations determined. The ${\rm gp}$-number of Cartesian products has been further studied in~\cite{klavzar-2019+}. In~\cite{klavzar-2019} the general posotion number has been connected with strong resolving graphs, and in~\cite{patkos-2019+} the general position number of Kneser graphs was investigated.

If $G$ is a graph and $\overline{G}$ its complement, then the {\em complementary prism} $G\overline{G}$ of $G$ is the graph formed from the disjoint union of $G$ and $\overline{G}$ by adding the edges of a perfect matching between the corresponding vertices of $G$ and $\overline{G}$~\cite{haynes-2007}. For example,  $C_5\overline{C}_5$ is the Petersen graph. Solely from this particular reason, but also from many additional ones, complementary prisms were studied from different perspectives. Since the Petersen graph is a key example in the theory of edge colorings, it is no surprise that the chromatic index of complementary prisms was studied in~\cite{zatesko-2019}. Other topics studied on complementary prisms include domination~\cite{haynes-2009}, cycle structure~\cite{meieling-2015}, complexity properties~\cite{duarte-2017}, spectral properties~\cite{cardoso-2018},  convexity number~\cite{castonguay-2019}, and b-chromatic number~\cite{bandeli-2019}. In this paper, we add to this list the general position problem. We proceed as follows. 

The next section contains definitions, observations, and known results needed in the rest of the paper. In Section~\ref{sec:upper} we prove that $\gp(\GG)$ is bounded from above by $n(G) + 1$ if $G$ is connected and by $n(G)$ if $G$ is disconnected, where $n(G)$ is the order of $G$. We also introduce the concept of $3$-general position sets and apply it to derive a characterization of graphs $G$ for which $\gp(\GG) = n(G) + 1$ holds provided that both $G$ and $\overline{G}$ are connected. Then, in Section~\ref{sec:lower}, we prove a sharp lower bound on $\gp(\GG)$. We follow with two sections on  complementary prisms of bipartite graphs and split graphs, respectively. In both cases, provided that a bipartite graph in question is connected, $\gp(G\overline{G})$ lies between $n(G)$ and $n(G)+1$. For connected bipartite graphs we characterize the graphs $G$ for which $\gp(G\overline{G})=n(G)+1$ holds, while for split graphs we give two partial results about the split graphs $G$  for which $\gp(G\overline{G})=n(G)$ holds. We conclude with Section~\ref{sec:block} in which we give a characterization of block graphs $G$ for which $\gp(G\overline{G})=n(G)+1$ holds and provide a family of block graphs in which the $\gp$-number of their complementary prisms is arbitrary smaller than their order.

\section{Preliminaries}
\label{sec:preliminary}

Graphs in this paper are finite and simple. Let $G = (V(G), E(G))$ be a graph. The maximum order of its complete subgraph is denoted by $\omega(G)$. Let further $\eta(G)$ denote the maximum order of an induced complete multipartite subgraph of the complement of $G$. The {\em distance} $d_G(u,v)$ between vertices $u$ and $v$ is the length of a shortest $u,v$-path. An $u,v$-path of minimum length is also called an $u,v$-{\it geodesic}. The {\em interval} $I_G[u,v]$ between $u$ and $v$  is the set of vertices that lie on some $u,v$-geodesic of $G$. For $S\subseteq V(G)$ we set $I_G[S]=\bigcup_{_{u,v\in S}}I_G[u,v]$. The {\it eccentricity} of $u$ is ${\rm ecc}_G(u) = \max \{d_G(u,v) : v\in V(G)\}$. The {\em radius} and the {\em diameter} of $G$ are ${\rm rad}(G) = \min \{{\rm ecc}_G(v):\ v\in V(G)\}$ and $\diam(G) = \max \{ecc_G(v):\ v\in V(G)\}$, respectively. A vertex $v$ is a \emph{central vertex} of $G$ if ${\rm ecc}_G(v)=\rad(G)$. The set of all central vertices is denoted by $C(G)$. We may simplify the above notation by omitting the index $G$ whenever $G$ is clear from the context. On the other hand, when we will want to emphasize that a vertex is central in a graph $G$, we will say that it is {\em $G$-central}. 

For a characterization of general position sets we needs some additional definitions. Let $G$ be a connected graph, $S\subseteq V(G)$, and ${\cal P} = \{S_1, \ldots, S_p\}$ a partition of $S$. Then ${\cal P}$ is \emph{distance-constant} if for any $i,j\in [p]$, $i\ne j$, the distance $d(u,v)$, where $u\in S_i$ and $v\in S_j$ is independent of the selection of $u$ and $v$.  If ${\cal P}$ is a distance-constant partition, and $i,j\in [p]$, $i\ne j$, then let $d(S_i, S_j)$ be the distance between a vertex from $S_i$ and a vertex from $S_j$. A distance-constant partition ${\cal P}$ is {\em in-transitive} if $d(S_i, S_k) \ne d(S_i, S_j) + d(S_j,S_k)$ holds for arbitrary pairwise different $i,j,k\in [p]$. Now all is ready to recall the announced characterization. 

\begin{theorem}
{\rm\cite[Theorem 3.1]{bijo-2019}}
\label{thm:gp-sets-characterization}
Let $G$ be a connected graph. Then $S\subseteq V(G)$ is a general position set if and only if the components of $G[S]$ are complete subgraphs, the vertices of which form an in-transitive, distance-constant partition of $S$.
\end{theorem}

We will also make use of the following two results. 

\begin{theorem}
{\rm \cite[Theorem 4.1]{bijo-2019}}
\label{thm:2.3}
If $\diam(G) = 2$, then $\gp(G) = \max\{\omega(G), \eta(G)\}$.
\end{theorem}

\begin{theorem} {\rm \cite[Theorem 5.1]{bijo-2019}}
\label{thm:2.2}  If $G$ is a connected, bipartite graph with $n(G)\ge 3$, then $\gp(G)\leq \alpha(G)$. Moreover, if $\diam(G) \in \{2,3\}$, then $\gp(G)=\alpha (G).$
\end{theorem}

Let $G$ be a graph and $G\overline{G}$ its complementary prism. Then we will consider $V(G\overline{G})$ as the disjoint union of $V(G)$ and $V(\overline{G})$. We will use the convention that if $u\in V(G)\cap V(G\overline{G})$, then its unique neighbour in $V(G\overline{G}) \cap V(\overline{G})$ will be denoted with $\overline{u}$ and called the {\em partner} of $u$ in $\overline{G}$. We will extend this notation to sets of vertices, that is, if $X\subseteq V(G)$, then the set of the partners of the vertices from $X$ will be denoted with $\overline{X}$. Since the complementation is an idempotent operation, $\overline{G}\, \overline{\overline{G}}$ is isomorphic to $G\overline{G}$. Note further that if $\diam(G) = 2$, then $\diam(\GG) = 2$, while if $G$ is an arbitrary connected graph, then $\diam(\GG) \le 3$. Note finally that $\GG$ is always connected, no matter whether $G$ is connected or not. 

\section{Upper bounds}
\label{sec:upper}

In this section we bound $\gp(\GG)$ from the above by $n(G) + 1$ for connected graphs $G$ and by $n(G)$ for disconnected graphs $G$, both bounds being sharp.  We also characterize the graphs $G$ for which $\gp(\GG) = n(G) + 1$ holds provided that both $G$ and $\overline{G}$ are connected. For this purpose, the concept of $3$-general position sets is introduced along the way. 

\begin{theorem}
\label{thm:upper} 
Let $G$ be a graph. 

(i) If $G$ is connected, then $\gp(G\overline{G})\leq n(G)+1$.

(ii)  If $G$ is disconnected, then $\gp(G\overline{G})\leq n(G)$.
\end{theorem}

\proof (i) Let $S\subseteq V(G\overline{G})$, where $|S| \ge n(G)+2$. By the pigeonhole principle there exist vertices $u,v \in V(G)$ such that $\{u, v, \overline{u},  \overline{v}\} \subseteq S$. Since either $uv\in E(G)$ or $\overline{u}\,\overline{v}\in E(\overline{G})$, we see that either $v \in I[u,\overline{v}]$ or $\overline{v} \in I[v,\overline{u}]$. It follows that $S$ is not a general position set and we can conclude that $\gp(G\overline{G})\leq n(G)+1$.

(ii) Let $G_1, \ldots, G_r$, $r\geq 2$, be the components of $G$. Assume that $S\subseteq V(G\overline{G})$, where $|S|\ge n(G)+1$, is a general position set of $G\overline{G}$. Using the pigeonhole principle again, there exists $v \in V(G)$ such that $\{v,\overline{v}\}\subseteq S$. We may without loss of generality assume that $v\in V(G_1)$. Let $x\in S\cap V(G)$, $x\ne v$. Then $xv\notin E(G)$, for otherwise $v$ would lie on a $x,\overline{v}$-geodesic. Moreover, $d_{\GG}(x,v) \le 2$, for otherwise, having in mind that $\overline{x}\,\overline{v}\in E(\GG)$, we would have that $\overline{v}$ would lie on a $x,v$-geodesic. It follows that $d_{\GG}(x,v) = 2$ holds. This in particular implies that $S\cap V(G) \subseteq V(G_1)$. Using a parallel argument we infer that if $\overline{y}\in S\cap V(\overline{G})$, then $d_{\GG}(\overline{y},\overline{v}) = 2$. This in turn implies that $\overline{y}\,\overline{v}\notin E(\overline{G})$ so that $yv\in E(G)$. It follows that $\{y:\ \overline{y}\in S\cap V(\overline{G})\} \subseteq N_G[v]$. 

We have thus proved that $S\subseteq V(G_1)\cup V(\overline{G}_1)$. Since $G_1$ is connected, the proof of (i) restricted to $G_1\overline{G}_1$ implies that $v$ is the unique vertex of $S\cap V(G)$ such that $\overline{v}\in S$. Hence $|S| \le n(G_1) + 1 \le n(G)$,  a contradiction. Hence $\gp(G\overline{G})\leq n(G)$.
\qed

To quickly demonstrate that the bounds of Theorem~\ref{thm:upper} are sharp, consider the following sporadic examples. First, as observed in~\cite{manuel-2018a} for the Petersen graph,  $\gp(C_5\overline{C_5}) = 6$, which demonstrates  sharpness of (i). Second, $\gp(\overline{P_2}\, \overline{{\overline{P_2}}}) = 2$  demonstrates sharpness of (ii). 

Our next goal is to characterize the graphs $G$ such that both $G$ and $\overline{G}$ are connected and $\gp(G\overline{G}) = n(G)+1$. By the above, the Petersen graph belongs to this family. 

\begin{lemma}
\label{lemma:3.3}
Let $G$ be a graph with $n(G)\ge 2$ and such that both $G$ and $\overline{G}$ are connected. If $\gp(G\overline{G})=n(G)+1$, then the following properties hold.
\begin{enumerate}
\item[(i)] $\rad(G)= 2$.
\item[(ii)] If $S$ is a $\gp$-set of $G\overline{G}$, then there exists a $G$-central vertex $v \in S$ such that $S\cap V(G)=\{u\in V(G):\ d_{G}(u,v)=2\}\cup \{v\}$ and $S\cap V(\overline{G}) = \overline{N_{G}[v]}$.
\end{enumerate}
\end{lemma}

\proof 
Let $S$ be a $\gp$-set of $G\overline{G}$. Then by lemma's assumption, $|S|=n(G)+1$, and hence there exists a vertex $v\in V(G)$ such that $\{v,\overline{v}\}\subseteq S$. 

 Suppose that there exists a vertex $x\in S\cap V(G)$ for which $d_G(x,v)\geq 3$ holds. Then the path on the vertices $x,\overline{x},\overline{v}, v$ is a $x,v$-geodesic in $G\overline{G}$ containing the vertex $\overline{v}$, a contradiction. Hence $d_G(x,v)\le 2$ for all $x\in S\cap V(G)$, $x\ne v$. Since $\overline{G}$ is connected, $V(G)\setminus N_G(v)\neq \emptyset$, for otherwise $\overline{v}$ would be an isolated vertex in $\overline{G}$. We have shown that ${\rm ecc}_G(v)= 2$ which in turn implies (i).

Theorem~\ref{thm:gp-sets-characterization} implies that $N_G(v)\cap S=\emptyset$ and $S\cap V(\overline{G})\subseteq \overline{N_G[v]}$.  Since $|S|=n(G)+1$ and because by the proof of Theorem~\ref{thm:upper}(i), the vertex $v$ is the unique vertex of $S\cap V(G)$ such that its partner in $\overline{G}$ also belongs to $S$, we must have that $S\cap V(G)=V(G)\backslash N_{G}(v)$ and $S\cap V(\overline{G})= \overline{N_{G}[v]}$. Combining this with the already proved property (i), the assertion (ii) follows. 
\qed

\begin{figure}[ht!]
\centering
\begin{pspicture}(0,0)(3,3)
 \psset{showpoints=true} 
\psline(0,0)(1,1)(0,2)
\psline(1,1)(2,1)(3,2) \psline(2,1)(3,0) \uput[d](0,0){$u_2$}
\uput[u](0,2){$u_1$} \uput[d](1,1){$u$} \uput[d](2,1){$v$}
\uput[u](3,2){$v_1$}\uput[d](3,0){$v_2$}
\end{pspicture}
\vspace{6mm} 
\caption{The double star $G$}
\label{fig1}
\end{figure}
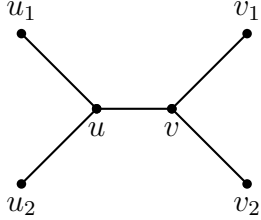

The condition (i) of Lemma~\ref{lemma:3.3} is not sufficient for $G$ to have $\gp(G\overline{G})=n(G)+1$. To see it, consider the double star $G$ shown in Fig.~\ref{fig1}. 

Clearly, $\rad(G)=2$ and both $G$ and $\overline{G}$ are connected. If $\gp (G\overline{G})=7$, and $S$ is a $\gp$-set of $G\overline{G}$, then by Lemma~\ref{lemma:3.3} either $v\in S$ or $u\in S$. Assume without loss of generality that $v \in S$. Using Lemma~\ref{lemma:3.3} again we get that $S = \{v, \overline{v}, u_1, u_2, \overline{u}, \overline{v_1}, \overline{v_2}\}$. But then the path on  the vertices $\overline{u}, u, v, \overline{v}$ is a geodesic containing three vertices of $S$, a contradiction. 

To state the announced characterization, we introduce the following concept. A set $S$ of vertices in a graph $G$ is a {\em $3$-general  position set} if no three vertices from $S$ lie on a common geodesic of length at most $3$. 

\begin{lemma}
\label{lem_3-general}
Let $G$ be a graph. A set $S\subseteq V(G)$ is a $3$-general position set if and only if the components of $G[S]$ are complete subgraphs, and if $d(u,v)=2$, where $u$ and $v$ lie in respective components $Q$ and $Q'$ of $G[S]$, then $d(x,y)=2$ for all $x\in Q$ and all $y\in Q'$.
\end{lemma}

\proof
Let $S\subseteq V(G)$ be a $3$-general position set. Then each component of $G[S]$ must be complete, for otherwise a non-compete component would contain an induced $P_3$ which is a geodesic. Suppose now that $d(u,v)=2$, where $u$ and $v$ lie in components $Q$ and $Q'$
of $G[S]$. Let $x\in Q$ and $y\in Q'$. Clearly, $1\le d(x,y)\le 3$. Suppose $d(x,y) = 1$. Assuming without loss of generality that $x\ne u$ we see that the vertices $y, x, u$ induce a geodesic, a contradiction.   Suppose $d(x,y) = 3$. Again assuming without loss of generality that $x\ne u$ we see that the vertices $x, u, v$ lie on a common geodesic of length $3$, another contradiction. Hence  $d(x,y)=2$.   

Conversely, suppose that the components of $G[S]$ are complete subgraphs and that the distance condition is fulfilled. Suppose on the contrary that $S$ is not a $3$-general position set, that is, there are three vertices $u, v, w\in S$ that lie on a common geodesic $P$ of length at most $3$. As the components of $G[S]$ are complete, $P$ cannot be of length $2$, so we may without loss of generality assume that $P$ in the path on vertices $u, v, x, w$ for some vertex $x\in V(G)$. Then $u$ and $w$ lie in different cliques of $G[S]$, say $u\in Q$ and $w\in Q'$. Clearly, then also $v\in Q$. Since $d(v,w) = 2$, the distance condition implies that $d(u,w) = 2$ as well. But this contradicts the fact that $P$ is a geodesic. 
\qed

We are now ready for the announced characterization. 

\begin{theorem}
\label{thm:3.5}
Let $G$ be a graph with $n(G)\ge 2$ and such that both $G$ and $\overline{G}$ are connected. Then $\gp(G\overline{G})=n(G)+1$ if and only if $\rad(G)= 2$ and there exists $v\in C(G)$ such that
\begin{enumerate}
\item[(i)] $\overline{N_{G}(v)}$ is a $3$-general position set in $\overline{G}$ and $N_{G}^{2}(v)$ is a $3$-general position set in $G$, and 
\item[(ii)] for each $x\in N_G(v)$ there exists $y\in N_{G}^{2}(v)$ such that $xy\notin E(G)$. 
\end{enumerate}
\end{theorem}

\proof 
First suppose that $\gp(G\overline{G})=n(G) + 1$. Let $S$ be a $\gp$-set of $G$, so that $|S| = n(G) + 1$. Then by Lemma \ref{lemma:3.3}, $\rad(G)=2$ and there exists a vertex $v\in C(G)$ such that $S\cap V(G)=N_{G}^{2}(v)$  and $S\cap V(\overline{G})=\overline{N_{G}(v)}$ are both general position sets in $G\overline{G}$. Hence $\overline{N_{G}(v)}$ is a $3$-general position set in $\overline{G}$ and $N_{G}^{2}(v)$ is a $3$-general position set in $G$. Thus (i) holds for the vertex $v$. Let now  $x\in N_G(v)$ and assume that $xy\in E(G)$ for every $y\in N_{G}^{2}(v)$. This implies that $d_{G\overline{G}}(\overline{x}, \overline{v}) = 3$. Since $d_{G\overline{G}}(\overline{x}, v) = 2$ and $v, \overline{v}$ are in the same clique of $G\overline{G}[S]$, Lemma~\ref{lem_3-general} implies that also $d_{G\overline{G}}(\overline{x}, \overline{v}) = 2$ must hold. This contradiction proves property (ii) for the vertex $v$. 

Conversely, suppose that $\rad(G)=2$ and that there exists a vertex $v$ in $C(G)$ satisfying properties (i) and (ii). We claim that the set $S=N_{G}^{2}(v)\cup \overline{N_G(v)}\cup \{v,\overline{v}\}$ induces an in-transitive, distance-constant partition into cliques. By Lemma~\ref{lem_3-general}, each of $N_{G}^{2}(v)$ and $\overline{N_G(v)}$ induces an in-transitive, distance-constant partition into cliques. Now, for any $u\in V(G)$ and $\overline{w}\in V(\overline{G})$ with $u\neq w$, we have $d_{G\overline{G}}(u,\overline{w})=2$. Which implies that $N_{G}^{2}(v)\cup \overline{N_G(v)}\cup\{v\}$ has an intransitive-distance constant partition into cliques. Next, for any $y\in N_{G}(v)$ we have $d_{G\overline{G}}(\overline{y},v)=2$ and by condition (ii) of our assumption, $d_{G\overline{G}}(\overline{y},\overline{v}) = 2$. Hence $S$ induces an in-transitive, distance-constant partition into cliques. Since ${\rm ecc}_G(v)=2$, we have that $|S|=n(G) + 1$.  Hence by Theorem~\ref{thm:gp-sets-characterization}, $S$ is a general position set of size $n(G) + 1$ and so $\gp(G\overline{G})=n(G)+1$ by Theorem~\ref{thm:upper}.  
\qed

As an example consider again the Petersen graph $P = C_5\, \overline{C_5}$. Let the vertices of the $C_5$ be  $v_1, v_2, v_3, v_4, v_5$. Observe that ${\rad}(C_5)=2$ and $C(C_5)=V(C_5)$. Select $v_1$ as a vertex from $C(C_5)$. Then $N_{C_5}^{2}(v_1)=\{v_3, v_4\}$ and $\overline{N_{C_5}(v_1)}=\{\overline{v_{2}},\overline{ v_{5}}\}$. It is now straightforward to check that the conditions of Theorem~\ref{thm:3.5} are fulfilled, hence $\gp(P) = 6$ with $\{v_1, \overline{v_1}, v_3, v_4, \overline{v_2}, \overline{v_5}\}$ being its $\gp$-set. 

Next, we give an infinite family of graphs that satisfies Theorem~\ref{thm:3.5}. Let $G$ be the graph obtained from the path on vertices $u, v, w$ and  disjoint cliques $K_{n_{1}}, \ldots, K_{n_{r}}$ and $K_{m_{1}}, \ldots, K_{m_{s}}$ by joining $u$ to all the vertices of $K_{n_{i}}$, $i\in [r]$, and $w$ to all the vertices of $K_{m_{j}}$, $j\in [s]$. Then $\rad(G)=2$ and $C(G)=\{v\}$. Since the sets $N_{G}(v)=\{u,w\}$ and $N_{G}^{2}(v)=V(G)\backslash \{u,v,w\}$ satisfy the  conditions of Theorem~\ref{thm:3.5}, we conclude that $\gp(G\overline{G})=n(G)+1$.

\section{A lower bound}
\label{sec:lower}

In this section we prove a sharp lower bound on $\gp(\GG)$. For this sake we define a new invariant $\overline{\gp}_{3}(G)$ of a graph $G$ as follows: 
$$\overline{\gp}_{3}(G) = \max \{ \gp_{3}(\overline{G}[V(\overline{G})\setminus \overline{S}]):\ S \gp_{3}{\mbox -}{\rm set\ of}\ G\}\,.$$

\begin{theorem}
\label{thm:lower}
If $G$ is a graph, then 
$$\gp(\GG) \ge \max \{\gp_{3}(G) + \overline{\gp}_{3}(G),  \gp_{3}(\overline{G}) + \overline{\gp}_{3}(\overline{G})\}\,.$$
Moreover, the bound is sharp. 
\end{theorem}

\proof
Let $S$ be a $\gp_{3}$-set of $G$ and let $Q_1, \ldots, Q_k$ be the complete subgraphs of $G[S]$ according to Lemma~\ref{lem_3-general}. For $x\in Q_i$ and $y\in Q_j$, where $i\ne j$, $d_{G}(x,y)\in \{2,3\}$ and hence $d_G(x,y) = d_{\GG}(x,y)$. This already implies that the cliques  $Q_1, \ldots, Q_k$ form an in-transitive, distant-constant partition of  $S$ in $\GG$. Consider now a $\gp_{3}$-set $\overline{T}$ of $\overline{G}[V(\overline{G})\setminus \overline{S}]$ and let $R_1, \ldots, R_t$ be the complete subgraphs of $\overline{G}[\overline{T}]$ according to Lemma~\ref{lem_3-general}. Then, by the same argument as we used for  $Q_1, \ldots, Q_k$, we infer that  the cliques  $R_1, \ldots, R_t$ form an in-transitive, distant-constant partition of  $\overline{T}$ in $\GG$.  Let now $x\in Q_i$, $i\in [k]$, and $y\in R_j$, $j\in [t]$. Then by the structure of $\GG$ it follows that $d_{\GG}(x,y) = 2$. If follows that  the cliques  $Q_1, \ldots, Q_k,  R_1, \ldots, R_t$ form an in-transitive, distant-constant partition in $\GG$.   With Theorem~\ref{thm:gp-sets-characterization} in hands we have thus proved that $\gp(\GG) \ge \gp_{3}(G) + \overline{\gp}_{3}(G)$. Starting with a $\gp_{3}$-set of $\overline{G}$ and repeating the argument from the paragraph above, we also get that $\gp(\GG) \ge \gp_{3}(\overline{G}) + \overline{\gp}_{3}(\overline{G})$.

By a simple argument we can see that if $n\ge 2$, then $\gp(K_n\overline{K}_n) = n$. Since $\gp_{3}(K_n) = n$ (as well as $\gp_{3}(\overline{K}_n) = n$), the bound is sharp. 
\qed

From Theorem~\ref{thm:lower} it is clear that $\gp(G\overline{G})\geq \max\{\gp_{3}(G),\gp_{3}(\overline{G})\}$. We next characterize the connected graphs for which the equality holds.

\begin{theorem}
	\label{thm:4.2}
	Let $G$ be a connected graph. Then  $\gp(G\overline{G})=\max\{\gp_{3}(G),\gp_{3}(\overline{G})\}$ if and only if $G$ is a complete multipartite graph.
\end{theorem}

\proof
Let  $\gp(G\overline{G})=\max\{\gp_{3}(G),\gp_{3}(\overline{G})\}$. 

Suppose first that  $\gp(G\overline{G})=\max\{\gp_{3}(G),\gp_{3}(\overline{G})\}=\gp_{3}(G)$. Let $S$ be a maximum 3-general position set of $G$. Then $S$ is also a maximum general position set of $G\overline{G}$. If $S\neq V(G)$, then $S\cup\{\overline{v}\}$ is a general position set of $G\overline{G}$ for all $v\in V(G)\backslash S$. This is impossible. Hence $\gp_{3}(G)=|S|=n(G)$. Then $G[S]$ must be connected and from Lemma~\ref{lem_3-general} we conclude that  $G=K_{n}$. 

Suppose second that $\gp(G\overline{G})=\max\{\gp_{3}(G),\gp_{3}(\overline{G})\}=\gp_{3}(\overline{G})$. If $\overline{G}$ is connected, then as above we have that $G=\overline{K}_n$. This is a contradiction to the fact that $G$ is connected. Hence $\overline{G}$ must be disconnected. Let $S$ be a maximum 3-general position set of $\overline{G}$. Then $S$ is also a maximum general position set of $G\overline{G}$.  If $S\neq V(\overline{G})$, then $S\cup\{v\}$ is a general position set of $G\overline{G}$ for all $\overline{v}\in V(\overline{G})\backslash S$. This is impossible and so $\gp_{3}(\overline{G})=|S|=n(\overline{G})$. By Lemma \ref{lem_3-general},  each component of $\overline{G}$ must be a clique. This shows that $G$ is a complete multipartite graph. 

Conversely, Let $G$ be a complete multipartite graph. Then $\overline{G}$ is a disjoint union of cliques. This shows that $V(G)$ is a general position set of $G\overline{G}$.  Hence it follows from Theorem~\ref{thm:upper}(ii) that $\gp(G\overline{G})=n(G)=\max\{\gp_{3}(G),\gp_{3}(\overline{G})\}$.
\qed

Note that if $\gp(G)\ge n(G)-2$, then Theorem~\ref{thm:lower} readily implies that $\gp(\GG) \ge n(G)$. On the other hand, the bound of Theorem~\ref{thm:lower} is never larger than $n(G)$, hence the existence of graphs $G$ for which the equality $\gp(\GG) = n(G) + 1$ holds (cf.\ Theorem~\ref{thm:3.5}), implies that the bound is not sharp in general. Additional sharpness cases for the bound of Theorem~\ref{thm:lower} will be presented in the subsequent sections. 

\section{Bipartite graphs}
\label{sec:bipartite}

In this section we give our attention to the complementary prisms of bipartite graphs. If $G=(V(G),E(G))$ is a bipartite graph with bipartition $V(G)=A\cup B$, then we write $G$ as a triple $(A,B,E(G))$. In $G=(A,B,E(G))$, set $$U_{G}=\{u\in A:\ \deg(u)=|B|\} \cup \{v\in B:\ \deg(v)=|A|\}\,.$$

\begin{theorem}
\label{thm:3.7} 
If $G=(A,B,E(G))$ is a connected, bipartite graph, then $n(G)\leq \gp(G\overline{G})\leq n(G)+1$. Moreover, $\gp(G\overline{G})=n(G)+1$ if and only if $\rad(G)=2$ and $C(G)$ is an independent set.
\end{theorem}

\proof 
By Theorem~\ref{thm:upper}, $\gp(G\overline{G})\leq n(G)+1$. Since independent sets and cliques are $3$-general position sets in any graph, it follows from Theorem~\ref{thm:lower} that the set $S=A\cup \overline{B}$ is a general position set in $G\overline{G}$. Thus $\gp(G\overline{G})\geq n(G)$.  It thus remains to characterize the graphs $G$ for which $\gp(G\overline{G})=n(G)+1$ holds. 

Suppose that $\gp(G\overline{G})=n(G)+1$. Then, in view of Theorem~\ref{thm:upper}, both $G$ and $\overline{G}$ are connected. Let $S$ be a general position set in $G\overline{G}$ of size $n(G)+1$. Then by Lemma \ref{lemma:3.3}, we have that $\rad(G)=2$ and there exists a $G$-central vertex $v\in S$ such that $S\cap V(G)=\{u\in V(G):\ d_{G}(u,v)=2\}$ and $S\cap V(\overline{G} ) = \overline{N_{G}[v]}$. Now, without loss of generality  we may assume that $v\in A$. Since ${\rm ecc}_{G}(v)=2$ and $G$ is bipartite, it follows that $v\in U_{G}$ and so $S\cap V(\overline{G})=\overline{N_{G}[v]}=\overline{B}\cup\{\overline{v}\}$. Now, we claim that $U_{G} \subseteq A$. Assume on the contrary that there exists a vertex $u\in B\cap U_{G}$. Then $d_{G\overline{G}}(\overline{u},\overline{v})=3$ and $\overline{u},\overline{v}\in S$. Moreover, $\overline{u},u,v,\overline{v}$ is a $\overline{u},\overline{v}$-geodesic in $G$ containing the vertex $v$. This leads to a contradiction to the fact that $S$ is a general position set in $G\overline{G}$. Thus $U_{G}\subseteq A$. This shows that $\rad(G)=2$ and that $C(G)$ is an independent set in $G$.

Conversely, suppose that $ \rad(G)=2$ and $C(G)$ is independent in $G$. Since $\rad(G)=2$ and $G$ is bipartite,  $C(G)=U_{G}$. Because $C(G)$ is independent, and $C(G)=U_{G}$, it follows that either $C(G)\subseteq A$ or $C(G)\subseteq B$, say $C(G)\subseteq A$. Let $v$ be a vertex in $C(G)$. We claim that the set $S=A\cup \overline{B}\cup \{\overline{v}\}$ is a general position set in $G\overline{G}$. Now, the set $N_{G}^{2}(v)=A\backslash \{v\}$ is an independent set in $G$ and so it is a $3$-general position set in $G\overline{G}$. Also, $\overline{N_{G}(v)}=\overline{B}$ is a clique in $\overline{G}$ and so it is a $3$-general position set in $G\overline{G}$. Moreover, since $v \in U_{G}=C(G)\subseteq A$, we have that $d_{G\overline{G}}(\overline{y},\overline{v})=2$ for all $y\in N_{G}(v)$. Hence it follows from Theorem~\ref{thm:3.5} that $S$ is a general position set in $G\overline{G}$ and hence $\gp(G\overline{G})=n(G)+1$.   
 \qed

In the rest of the section we present the general position number of complementary prism of some standard families of bipartite graphs. Let $T$ be a tree. Then $\gp(T)$ is the number of its leaves~\cite{bijo-2019, ullas-2016}, and $\gp(\overline{T})=\max\{\alpha(T),\triangle(T)+1\}$~\cite{bijo-2019}. Since $1\le |C(T)|\le 2$, Theorem~\ref{thm:3.7} implies that $\gp(T\overline{T})=n(G)+1$ if and only if $|{\rm C}(T)|=1$ and $\rad(T)=2$. This is possible only when $T$ has diameter $4$. We have thus deduced: 

\begin{corollary}
\label{corollary:3.8}
If  $T$ is a tree, then  
$$\gp(T\overline{T}) = 
\begin{cases}
 n(G)+1; & \diam(T) = 4\,, \\
 n(G); & otherwise\,.
 \end{cases}
$$
\end{corollary}

The {\em Cartesian product} $G \cp H$ of graphs (factors) $G$ and $H$ has $V(G \cp H) = V(G)\times V(H)$ and  vertices $(g,h)$ and $(g',h')$ are adjacent if either $g=g'$ and $hh'\in E(G)$, or $h=h'$ and $gg'\in E(G)$. For $n,m\geq 2$ set $P_{nm} = P_n \square P_m$. In~\cite{manuel-2018b} it was proved that $\gp(P_{nm})=4$ for $n,m\geq 3$, while in~\cite{bijo-2019} the following result was deduced: 
$$\gp(\overline{P}_{nm})=
\begin{cases}
 4; & n = m = 2\,,\\
 \lceil \frac{n}{2}\rceil\lceil\frac{m}{2}\rceil + \lfloor\frac{n}{2}\rfloor\lfloor\frac{m}{2}\rfloor; & otherwise\,.
 \end{cases}$$
With the help of Theorem~\ref{thm:3.7}, we can add to these results the following. 

\begin{corollary}
\label{corollary:3.9}
If $n,m\geq 2$,  then
$$\gp(P_{nm}\overline{P}_{nm})=
\begin{cases}
  10; & n = m = 3\,, \\
   nm; & otherwise\,.
 \end{cases}$$
 \end{corollary}

As the last subclass of bipartite graphs consider hypercubes. Recall that the $n$-cube $Q_n$ is the $n$-fold Cartesian product of $K_2$.  Once more applying Theorem~\ref{thm:3.7} we get: 

 \begin{corollary}
\label{corollary:3.10}
If $n\geq 2$, then $\gp(Q_n\overline{Q}_n) = 2^n$.
\end{corollary} 

\section{Split graphs}
\label{sec:split}

A graph $G=(V(G),E(G)) $ is a {\em split graph} if $V(G)$ can be partitioned into a clique $C$ and an independent set $I$. If so, the pair $(C,I)$ is a {\em split partition} of $G$ and we write $G=(C,I,E(G))$.

\begin{theorem}
\label{thm:3.12} 
If $G=(C,I,E(G))$ is a split graph, then $n(G)\leq gp(G\overline{G})\leq n(G)+1$. Moreover, the following  hold. 
 \begin{enumerate}
  \item[(i)] If $\deg_{G}(x)\geq |C|+1$ for all $x\in C$ and $\deg_{G}(y)\leq |C|-2$ for all $y\in I$, then $\gp(G\overline{G})=n(G).$
  \item[(ii)]If $\gp(G\overline{G})=n(G)$, then $\overline{G}$ is disconnected or $\deg_{G}(x)\geq |C|$ for all $x\in C$.
  \end{enumerate} 
 \end{theorem}

\proof
By Theorem~\ref{thm:lower}, the set $H=C\cup \overline{I}$ is a general position set of $G\overline{G}$ and hence $gp(G\overline{G})\geq n(G)$. The upper bound again follows from Theorem~\ref{thm:upper}. 
 
 (i) Suppose that $\deg_{G}(x)\geq |C|+1$ for all $x\in C$ and $\deg_{G}(y)\leq |C|-2$ for all $y\in I$. By way of contradiction suppose that $\gp(G\overline{G})=n(G)+1$. Then by Lemma~\ref{lemma:3.3}, there exists a $G$-central vertex $v$ such that $S=N_{G}^{2}(v)\cup \overline{N_G(v)}\cup \{v,\overline{v}\}$ is a general position set of $G\overline{G}$. We consider two cases. 
 
\medskip\noindent
{\bf Case 1}: $v\in C$. \\
Then $\deg_{G}(v)\geq |C|+1$ and let $z$ be a neighbour of $v$ in $I$. Now, since $\overline{N_{G}[v]}\subseteq S$, we have that $\overline{z}\in S$. By our hypothesis, $\deg_{G}(z)\leq |C|-2$. Thus we can choose two distinct vertices, say $u_{1}$ and $u_{2}$ in $C$, such that both $u_{1}$ and $u_{2}$ are non-adjacent to $z$ in $G$. This shows that $\overline{z}\in I_{G\overline{G}}[\overline{u}_{1},\overline{u}_{2}]\subseteq I_{G\overline{G}}[S]$. This is a contradiction to the fact that $S$ is a general position set in $G\overline{G}$.

\medskip\noindent
{\bf Case 2}: $v\in I$. \\
In this case $\deg_{G}(v)\leq |C|-2$. Let $u$ be a vertex in $C$ such that $u$ and $v$ are non adjacent in $G$. Then by Lemma~\ref{lemma:3.3}, $u\in S$. Then $\deg_{G}(u)\geq |C|+1$, thus we can choose two distinct vertices, say $v_{1}$ and $v_{2}$ in $I$ such that both $v_{1}$ and $v_{2}$ are adjacent to $u$ in $G$. This shows that $u\in I_{G\overline{G}}[v_{1},v_{2}]\subseteq I_{G\overline{G}}[S]$. This is a contradiction to the fact that $S$ is a general position set in $G\overline{G}$. 

\medskip
Since in both cases we got a contradiction, we conclude that $\gp(G\overline{G})=n(G)$.
  
\medskip
(ii) Assume that $\gp(G\overline{G})=n(G)$ and $\overline{G}$ is connected. Then $G$ has no universal vertex and so $\rad(G)=2$. Suppose that $\deg_{G}(v)=|C|-1$ for some $v\in C$. Then $v$ is a $G$-central vertex and $\overline{N_{G}(v)}=\overline{C}$ is an 3-general position set in $\overline{G}$ and $N_{G}^{2}(v)=I $ is a 3-general position set of $G$. Moreover, since $G$ has no universal vertices, we have that for each $x\in C$ there exists $y\in I$ such that $xy \notin E(G)$. Hence by Theorem~\ref{thm:3.5}, $\gp(G\overline{G})=n(G)+1$, a contradiction. 
\qed
 
\begin{figure}[ht!]
\begin{minipage}{0.48\textwidth}

\centering
\begin{pspicture}(0,0)(3,3)
 \psset{showpoints=true} \psline(0,0)(1,1)(0,2) 
  \psline(0,2)(3,2) \psline(0,0)(3,0) \psline (0,2)(0,0)  \psline(1,1)(2,1)
\uput[u](0,2){$u_2$} \uput[u](1,1){$u_1$} \uput[d](2,1){$v_1$} \uput[d](0,0){$u_3$}
\uput[u](3,2){$v_2$}\uput[d](3,0){$v_3$}
\end{pspicture}
${G_{1}}$
\end{minipage}\hfill
\begin{minipage}{0.48\textwidth}
\centering
\begin{pspicture}(0,0)(3,3)
 \psset{showpoints=true} \psline(0,0)(1,1)(0,2) 
 \psline(2,1)(0,2) \psline(0,2)(3,2) \psline(0,0)(3,0) \psline (0,2)(0,0)  \psline(1,1)(2,1)
\uput[u](0,2){$u_2$} \uput[u](1,1){$u_1$} \uput[d](2,1){$v_1$} \uput[d](0,0){$u_3$}
\uput[u](3,2){$v_2$}\uput[d](3,0){$v_3$}
\end{pspicture}
${G_{2}}$
\end{minipage}
\vspace{6mm} 
\hspace{0.5cm}
\caption{Split graphs}
\label{fig:two-split-graphs}
\end{figure}
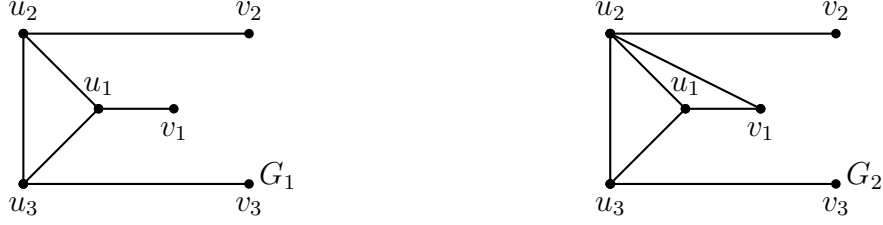

The converse for neither of Theorem~\ref{thm:3.12}(i) and (ii) is true. For this sake consider the split graphs $G_1$ and $G_2$ shown in Fig.~\ref{fig:two-split-graphs}.

Consider first the split graph $G_{1}$ with $C = \{u_1,u_2,u_3\}$ anf $I = \{v_1, v_2, v_3\}$. The center of $G_1$ is  $C(G_{1})=\{u_{1},u_{2},u_{3}\}$.  Since no central vertex satisfies conditions of Theorem~\ref{thm:3.5}, $\gp(G_{1}\overline{G}_{1}) = n(G_{1}) = 6$. Also, since $\deg_{G_{1}}(x)=3<4=|C|+1$ for all $x\in C(G_{1})$, we infer that the  converse of Theorem~\ref{thm:3.12}(i) does not hold. Next, in $G_{2}$, the vertex $u_{1}$ is a $G_{2}$-central vertex which satisfies the conditions of Theorem~\ref{thm:3.5}. Hence $\gp(G_{2}\overline{G}_{2})=7=n(G)+1$, whereas both $G_{2}$ and $\overline{G}_{2}$ are connected and have no simplicial vertices in the corresponding cliques. Hence the  converse of Theorem~\ref{thm:3.12}(ii) also does not hold. 

Again, consider a split graph $G=(C,I,E(G))$ such that both $G$ and $\overline{G}$ are connected. Then it is clear from Theorem~\ref{thm:3.12} that $\gp(G\overline{G})=n(G)+1$ if either $C$ or $\overline{I}$ contains a simplicial vertex.

We conclude the section with the natural problem arising from Theorem~\ref{thm:3.12}. 

\begin{problem}
Characterize the split graphs $G$ for which $\gp(G\overline{G})= n(G)+1$ holds. 
\end{problem}

\section{Block graphs}
\label{sec:block}

A graph is a {\em block graph} if every maximal 2-connected component is a clique. In a block graph each vertex is either a simplicial vertex or a cut vertex. The following theorem gives a characterization of block graphs for which the bound in Theorem~\ref{thm:upper}(i) is attained.

\begin{theorem}
\label{thm:3.14}
Let $G$ be a block graph. Then $\gp(G\overline{G})=n(G)+1$ if and only if $\rad(G)=2$ and there exists a $G$-central vertex $v$ such that $N_{G}(v)$ is either a clique or an independent set, and $N_{G}(v)$ contains at least two cut vertices of $G$.
\end{theorem}

\proof 
First suppose that $\gp(G\overline{G})=n(G)+1$. Then by Lemma~\ref{lemma:3.3}, we have $\rad(G)=2$ and $G$ contains a central vertex $v$ such that $N_{G}^{2}(v)\cup \overline{N_{G}(v)}\cup\{v, \overline{v}\}$ is a $\gp$-set of $G\overline{G}$. We claim that $N_{G}(v)$ is either a clique or an independent set. Let $B_{1}, \ldots, B_{k}$ be the  blocks of $G$ containing $v$. Then we have $N_{G}(v)=\cup^{k}_{i=1}N_{B_{i}}(v)$. Since each $N_{B_{i}}(v)$ induces  a clique, $\overline{N_{G}(v)}$ induces a multipartite set.  By Theorem~\ref{thm:3.5}, we have that $\overline{N_{G}(v)}$ is a 3-general position set in $\overline{G}$. Thus either $k=1$ or $\overline{N_{G}(v)}$ is a clique. This proves the claim. Next we  prove that $N_{G}(v)$ contains at least two cut vertices of $G$. Since $v$ is a central vertex and $\rad(G)=2$, we have that $|N_{G}(v)|\geq 2$. Let $y\in N_{G}(v)$. Since $\ecc_G(v)=2$, it follows from Theorem~\ref{thm:3.5} that corresponding to the vertex $y$ there exist vertices $x\in N_{G}^{2}(v)$ and $w\in N_{G}(v)$ such that $xw\in E(G)$ and $xy\notin E(G)$. Similarly, corresponding to the vertex $w$ there exist vertices $x'\in N_{G}^{2}(v)$ and $z\in N_{G}(v)$ such that $x'z\in E(G)$ and  $x'w\notin E(G)$. Since $G$ is block graph, both $w$ and $z$ must be cut vertices in $G$.

Conversely, suppose that $\rad(G)=2$ and $G$ contains a central vertex $v$ such that $N_{G}(v)$ is either a clique or an independent set, and $N_{G}(v)$ contains at least two cut vertices of $G$. Since $\overline{N_{G}(v)}$ is either a clique or an independent set in $G$, it is a 3-general position set of $\overline{G}$. Let $B'_{1}\ldots ,B'_{r}$ be the blocks of $G$ which do not contain $v$. Then each  $B'_{i}$, $i\in [r]$, must contains a vertex from $N_G(v)$, say $u_{i}$. Then $N_{G}^{2}(v) = \cup^{r}_{i=1}V(B'_{i})\setminus \{u_i\}$. Since $V(B'_{i}) \setminus \{u_i\}$ induces a clique in $G$, we have that $N_{G}^{2}(v)$ is a union of disjoint cliques in $G$. This implies that $N_{G}^{2}(v)$ is a 3-general position set in $G$. We now claim that for each $y\in N_{G}(v)$, there exists a vertex $x$ in $N_{G}^{2}(v)$ such that $xy \notin E(G)$. Suppose on the contrary a vertex $y\in N_{G}(v)$ is adjacent to each vertex in $N_{G}^{2}(v)$. But this would mean that $N_G(v)$ has at most one cut vertex, a contradiction. Hence, by Theorem~\ref{thm:3.5}, we   conclude that  $N_{G}^{2}(v)\cup \overline{N_{G}(v)}\cup\{v, \overline{v}\}$ is a $\gp$-set of $G\overline{G}$ and so $\gp(G\overline{G})=n(G)+1$.
\qed

In the rest of the section we present an infinite family of block graphs such that the $\gp$-number of their complementary prisms is arbitrary smaller than their order. Let $G_{0}$ be the complete graph $K_{3}$ with the vertex set $\{v_{1},u_{1},v_{2}\}$. For $k\geq1$, the graph $G_{k}$ is obtained from $G_{k-1}$ by adding two new adjacent vertices $u_{k+1}$ and $v_{k+2}$ and joining both $u_{k+1}$ and $v_{k+2}$ to the vertex $v_{k+1}$. Note that $G_{k}$ is a block graph with $k+1$ blocks $B_{1},\ldots ,B_{k+1}$, where each block $B_i$ is a triangle with $V(B_{i})= \{v_{i},u_{i},v_{i+1}\}$. Note further that $E_{k} = \{v_1, v_{k+2} \} \cup \{u_{1}, \ldots, u_{k+1}\}$ is the set of simplicial vertices of $G_k$ and that the set of remaining vertices, that is $A_k = \{v_{2}, \ldots, ,v_{k+1}\}$, is the set of cut vertices of $G_{k}$. Clearly, $A_{k}$ induces a path of length $k-1$, let $X_k$ and $Y_k$ form the bipartition of $A_{k}$ where $|X_k| \geq |Y_k|$.
 
\begin{theorem} 
If $k\geq 5$, then $\gp(G_{k}\overline{G}_{k})=n(G_k)-\lfloor \frac{k}{2} \rfloor$.
\end{theorem}

\proof
Let $k\ge 5$ and let $E_k$, $A_k$, $X_k$, and $Y_k$ be the sets of vertices as defined before the theorem. Using the proof of Theorem~\ref{thm:lower} we infer that $S=E_{k}\cup \overline{X}_k$ is a general position set of $G_{k}\overline{G}_{k}$ of size $n(G_k)-\lfloor \frac{k}{2} \rfloor $. Hence $\gp(G_{k}\overline{G}_{k})\geq n(G_k)- \lfloor \frac{k}{2} \rfloor $. 

To prove the other inequality, assume on the contrary that there exists a general position set $S$ of $G_{k}\overline{G}_{k}$ with $|S|> n(G_k)-\lfloor \frac{k}{2} \rfloor $. Let $M= S\cap V(G_{k})$ and $\overline{N}=S\cap V(\overline{G}_{k})$. Then both $M$ and $\overline{N}$ are $3$-general position sets of $G_{k}$ and $\overline{G}_{k}$, respectively. Since $|S|>n(G_k)-\lfloor \frac{k}{2} \rfloor$, it follows that either $|M|>k+3$ or $|N|=|\overline{N}|>\lfloor \frac{k}{2} \rfloor$.

\medskip\noindent
{\bf Claim A}:  $\gp_3(G_{k})=k+3$ for all $k\geq 0$. \\
We proceed by induction on $k$, the cases $k=0$ and $k=1$ being easily verified. Assume that $\gp_3(G_{i})=i+3$ for all $i$ with $2\leq i < k$ and consider $G_{k}$. Let $H_{k}$ be a maximum $3$-general position set of $G_{k}$ and suppose that $|H_{k}|>k+3$. Then both $u_{k+1},v_{k+2} \in H_{k}$. Otherwise, if $u_{k+1}\notin H_{k}$ (say), then $H_{k}\backslash \{v_{k+2}\}$ is a $3$-general position set of $G_{k-1}$ of size at least $k+3$. This is a contradiction to the induction hypothesis. Now, let $H_{k-1}= H_{k}\backslash \{u_{k+1},v_{k+2}\}$. Then $H_{k-1}$ is a $3$-general position set of $G_{k-1}$ of size at least $k+2$. hence by induction $H_{k-1}$ is a maximum $3$-general position set of $G_{k-1}$ of size $k+2$. Again, since $\gp(G_{k-2})=k+1$, it follows that either $u_{k}\in H_{k-1}$ or $v_{k+1}\in H_{k-1}$. Recall that $H_{k}$ is a $3$-general position set in $G_{k}$. This shows that $H_{k-1}$ contains exactly one vertex from the set $\{u_{k},v_{k+1}\}$, say $u_{k}$. But this leads to the fact that $H_{k-1}\cup \{v_{k+1}\}$ is a $3$-general position set of $G_{k-1}$ of size $k+3$, a contradiction, and the claim is proved. 

\medskip
Since $M$ is a $3$-general position set in $G_{k}$, Claim A implies that $|M|\leq k+3$. This shows that $|N|>\lfloor{ \frac{k}{2}} \rfloor$. We consider the following three cases. 

\medskip\noindent
{\bf Case 1}: $N\subseteq E_{k}$.\\
Recall that $N$ contains at most one vertex from $M$. First suppose that $M\cap N \neq \emptyset $. In this case there exists a simplicial vertex $x$ such that $x\in M$ and $\overline{x}\in \overline{N}$. Let $y$ be a vertex from $N$ distinct from $x$. If $x$ and $y$ are non-adjacent in $G_{k}$, then $x,\overline{x},\overline{y}$ is a shortest path in $G_{k}\overline{G}_{k}$ with $x,\overline{x},\overline{y}\in M\cup \overline{N}=S$. This is impossible. Hence $x$ and $y$ are adjacent in $G_{k}$. Now, since $x,y\in N\subseteq E_{k}$, it follows that either $x,y \in B_{1}$, or $x,y\in B_{k+1}$. This shows that $\overline{E_{k}\backslash \{x,y\}}\subseteq I_{G_{k}\overline{G}_{k}}[\overline{x},\overline{y}]\subseteq I_{G_{k}\overline{G}_{k}}[\overline{N}]$  and so $N=\{x,y\}$. This is a contradiction to the fact that $|N|>\lceil \frac{k}{2} \rceil \geq 3$. So, assume that $M\cap N=\emptyset $. Because $\overline{N}$ is a $3$-general position set in $\overline{G}_{k}$, we have that the components of the induced subgraph of $\overline{N}$ in $\overline{G}_{k}$ are cliques. Hence $N$ induces a complete multipartite graph in $G_{k}$. Since $N\subseteq E_{k}$, it follows that $N\subseteq E_{k}\backslash \{v_{1},v_{k+1}\}$ and so $|N|\leq k+1$. On the other hand, since $M$ is a $3$-general position set in $G_{k}$ and $A_{k}$ induces a path of order $k$, it follows that $|M\cap A_{k}|\leq \lceil \frac{k}{2} \rceil.$
 
Thus $|S|= |M|+|\overline{N}|=|M|+|N|=|M\cup N|=|(M\cup N)\cap E_{k}|+|(M\cup N)\cap A_{k}|\leq k+\lceil \frac{k}{2}\rceil +3=n(G_{k})-\lfloor \frac{k}{2}\rfloor$.

\medskip\noindent
{\bf Case 2}:  $N\subseteq A_{k}$. \\
Recall that $N$ induces a complete multipartite graph in $G_{k}$ and $A_{k}$ induces a path of order $k$. Hence it is clear that either $N\subseteq X_k$ or $N\subseteq Y_k$ and so $|N|\leq \lceil \frac{k}{2} \rceil$. Hence this case cannot occur.

\medskip\noindent
{\bf Case 3}: $N\cap E_{k}\neq \emptyset$ and $N\cap A_{k} \neq \emptyset$.

\medskip\noindent
\textbf{Claim B}: $N$ is an independent set in $G_{k}$ for all $k\geq 5$. \\
Assume on the contrary that there exists adjacent vertices $u$ and $v$ in $N$ and consider the following cases. 

\medskip\noindent
{\bf Subcase 3.1}: $u$ is a simplicial vertex and $v$ is a cut vertex. \\
We may assume that $u=u_{i}$ and $v=v_{i+1}$. Then $I_{G_{k}\overline{G}_{k}}[\overline{u},\overline{v}]$ covers both the $\overline{E_{k}\backslash\{u_{i+1}\}}$ and $\overline{A_{k}\backslash \{v_{i},v_{i+2}\}}$. This shows that $N\subseteq \{v_{i},v_{i+1},v_{i+2},u_{i},u_{i+1}\}$. Again, since $N$ induces a complete multipartite graph, it follows that $N=\{v_{i},v_{i+1},u_{i}\}$ and so $k=4$.

\medskip\noindent
{\bf Subcase 3.2}: both $u$ and $v$ are simplicial vertices. \\
In this situation either $u,v\in B_{1}$ or $u,v\in B_{k+1}$, say $u,v\in B_{1}$. Then $u=v_{1}$ and $v=u_{1}$. As in Subcase 3.1 we can then prove that $N\subseteq \{v_{1},v_{2},u_{1}\}$ and so $k\leq 4$.

\medskip\noindent
{\bf Subcase 3.3}: both $u$ and $v$ are cut vertices. \\
Then $u=v_{i}$ and $v=v_{i+1}$ for some $i\geq 2$. Similarly as in Subcase 3.1 we have that $N\subseteq \{v_{i-2},v_{i-1},v_{i},v_{i+1},u_{i},u_{i+1}\}$. Again, since $N$ induces a complete multipartite graph, it follows that $|N|\leq 3$ and so $k\leq 4$. Hence Claim B follows.

\medskip
Since $N$ is an independent set in $G_{k}$, we have that $N$ contains at most one vertex from each block and so $|N|\leq k+1$. Recall that $N$ contains both cut vertices and simplicial vertices. This shows that $|N|\leq k$. Again, since $N$ is an independent set, it follows that $M\cap N=\emptyset$. 

\medskip\noindent
{\bf Claim C}: If $v_{2}\in N$ or $v_{k+1}\in N$, then $|N|\leq \lceil {\frac{k}{2}}\rceil$ for all $k\geq 2$. \\
To prove this claim, we use induction on the number of cut vertices $k$. If $k=2$ or 3 the result holds. Assume the result holds for all integers $i$ with $3\leq i < k$ and consider $G_{k}$. Suppose that $v_{2}\in N$ or $v_{k+1}\in N$, say $v_{k+1}\in N$. Now, since $N$ is independent in $G_{k}$, we have that both $u_{k+1} \notin N$ and $v_{k+2}\notin N$. Hence $N\subseteq V(G_{k-1})$ and so by induction hypothesis $|N|\leq\lceil{\frac{k-1}{2}}\rceil \leq \lceil\frac{k}{2}\rceil$. Hence Claim C follows.

In the following, we prove that $|N| \leq \lceil \frac{k}{2}\rceil $. Choose $\alpha$ with $2\leq \alpha \leq k+1$ such that the cut vertex $v_{\alpha} \in N$. Then it follows from Claim C that  
$$| N\cap V(G_{\alpha-1})|\leq  \left\lceil {\frac{\alpha-1}{2} }\right\rceil \,.$$ and
$$|N\cap [(V(G_{k})\backslash V(G_{\alpha-3}))\cup \{v_{\alpha-1}\}]| \leq \left\lceil { \frac{k+1-\alpha}{2}} \right\rceil \,.$$ Since $v_{\alpha}$ is counted twice, we have that
\begin{align*}
| N| & \leq |N\cap V(G_{\alpha-1})|+| N \cap [(V(G_{k})\backslash V(G_{\alpha-3}))\cup \{v_{\alpha-1} \} ]| \\
& \leq  \left\lceil\frac{\alpha-1}{2}\right\rceil + \left\lceil {\frac{k+1-\alpha}{2}}\right\rceil -1\\
&\leq \left\lceil\frac{k}{2}\right\rceil\,.
\end{align*}
This shows that $|M\cup N|\leq k+\lceil \frac{k}{2} \rceil +3$.

Since in all the cases we have arrived at a contradiction, we conclude that $\gp(G_{k}\overline{G_{k}})=n(G_k)- \lfloor \frac{k}{2} \rfloor $. 
\qed

\section*{Acknowledgments}

Sandi Klav\v{z}ar acknowledges the financial support from the Slovenian Research Agency (research core funding P1-0297 and projects J1-9109, J1-1693, N1-0095). Neethu P K acknowledges the Council of Scientific and Industrial Research(CSIR), Govt. of India for providing financial assistance in the form of Junior Research Fellowship. Manoj Changat acknowledges the financial support from the Department of Science and Technology(research project under 'MATRICS' scheme No. MTR/2017/000238).

\end{document}